\documentclass[reqno]{amsart}

\usepackage{amsmath,amsthm,amssymb,tikz,hyperref,cleveref,CJKutf8}
\hypersetup{hidelinks}

\crefname{theorem}{Theorem}{Theorems}
\crefname{corollary}{Corollary}{Corollaries}
\crefname{figure}{Figure}{Figures}
\crefname{remark}{Remark}{Remarks}
\crefname{example}{Example}{Examples}
\crefname{section}{Section}{Sections}
\crefname{table}{Table}{Tables}

\definecolor{oi-sky}{HTML}{56B4E9}
\definecolor{oi-green}{HTML}{009E73}
\definecolor{oi-yellow}{HTML}{F0E442}
\definecolor{oi-blue}{HTML}{0072B2}
\definecolor{oi-vermillion}{HTML}{D55E00}

\usetikzlibrary{decorations.pathreplacing}
\tikzset{
  vtx/.style={circle, draw=black, thick, minimum size=6.5mm, inner sep=0pt,
    font=\small},
  deg3/.style={vtx, fill=oi-blue, text=white},
  deg2/.style={vtx, fill=oi-sky},
  deg1/.style={vtx, fill=oi-yellow},
  plain/.style={vtx, fill=white, minimum size=5mm},
  edge/.style={thick},
}

\newtheorem{theorem}{Theorem}
\newtheorem{corollary}[theorem]{Corollary}
\theoremstyle{definition}
\newtheorem{remark}[theorem]{Remark}
\newtheorem{example}[theorem]{Example}

\newcommand{\R}{\mathbb{R}}
\numberwithin{equation}{section}

\title{Degree-respecting labellings with ordered neighbour sums}
\author{Toshio Ishikawa}
\address{Independent Researcher}
\email{t.ishikawa.1017@gmail.com}

\author{Tomoki Uda}
\address{Department of Electronics and Communication Technology, Faculty of
Science and Technology, Nanzan University, Nagoya, Japan}
\email{uda0@nanzan-u.ac.jp}
\date{}


\begin{document}
\begin{CJK*}{UTF8}{min}

\begin{abstract}
In an exhaustive enumeration of hydrocarbons with only single bonds, carbon
skeletons are generated as adjacency matrices and pruned by two rules: a vertex
of larger degree receives a smaller label, and among vertices of equal degree,
a vertex whose adjacency row is the larger binary integer receives a smaller
label. The pruning is sound only if every skeleton admits a labelling obeying
both rules. We prove that such a labelling exists for every finite simple
graph, connected or not and with no bound on the degrees, and that the same
holds when the weights \(2^i\) of the labels \(i\) are replaced by any
strictly increasing real weights. The proof takes, among the labellings that
order vertices by non-increasing degree, one minimising the sum over edges of
the products of the endpoint weights.
\end{abstract}

\maketitle

\section{Introduction}

The question studied here arose in an exhaustive enumeration of hydrocarbons
with only single bonds, described by one of the authors (T.I.) in a blog post
\cite{smooth-pudding-2024}; the problem goes back to a post by T.I. dated
29 May 2024, cited there. Up to stereoisomerism,
such a molecule is determined by its carbon skeleton: the finite simple graph whose vertices are
the carbon atoms and whose edges are the carbon--carbon bonds. The skeleton is
connected, every vertex has degree at most \(4\), and hydrogen atoms fill the
remaining valences, so they need not be recorded. Enumerating these
hydrocarbons with \(n\) carbon atoms therefore amounts to listing, up to
isomorphism, the connected graphs on \(n\) vertices with maximum degree at
most \(4\).

The enumeration generates such graphs as adjacency matrices. We describe it
with the label set
\[
  I_n=\{0,1,\ldots,n-1\}
\]
in place of the labels \(1,\ldots,n\) used in \cite{smooth-pudding-2024}; by
\cref{rem:conventions}, this changes neither of the rules below. A bijective labelling of the vertices by \(I_n\) turns a graph into a symmetric
\(0/1\) matrix with zero diagonal, and different labellings of one graph
generally give different matrices. In practice, listing the adjacency matrices
of all such graphs under all labellings therefore produces each graph many
times, and every
duplicate has to be removed later by an isomorphism test. To reduce this
redundancy, the enumeration prunes while a matrix is being generated, as
implemented in \cite[directory \texttt{p2\_ch03-2}]{hydrocarbon-rta-2025}. It reads
the row of a vertex \(v\) as the \(n\)-bit integer whose bit in position \(j\),
of weight \(2^j\), is \(1\) exactly when \(v\) is adjacent to the vertex
labelled \(j\), and it keeps only matrices whose labelling obeys two rules.
\begin{enumerate}
  \item A vertex of larger degree receives a smaller label.
  \item Among vertices of equal degree, a vertex whose row is the larger
    integer receives a smaller label.
\end{enumerate}
For a labelling \(f\) and the set \(N(v)\) of neighbours of \(v\), the row
value of \(v\) is
\[
  \sum_{u\in N(v)}2^{f(u)},
\]
the sum of the weights \(2^i\) of the labels \(i\) carried by the neighbours of
\(v\); \cref{fig:motivation} shows these row values for a labelled
five-vertex skeleton that obeys both rules. The second rule thus requires that, among vertices of equal degree, the
vertex with the larger label have a row value no larger than that of the
other vertex.

The pruning is sound only if every skeleton admits at least one labelling
obeying both rules; otherwise that skeleton would be missing from the output. This existence had not been proved
when \cite{smooth-pudding-2024} was written. \Cref{cor:powers-of-two} below
proves it, in greater generality than the enumeration requires: the proof uses
neither connectedness nor the degree bound, and
\cref{thm:strict-weights} replaces the weights \(2^i\) by the values \(c(i)\)
of an arbitrary strictly increasing function \(c:I_n\to\R\).
\Cref{sec:remarks} discusses the lack of uniqueness, counts
for skeletons with \(2\) to \(7\) carbon atoms the adjacency matrices that remain
after each rule, and describes a check of the proofs in Lean~4.

\begin{figure}[t]
  \centering
  \begin{tikzpicture}[line cap=round, line join=round]
    \begin{scope}[yshift=-0.3cm]
      \node[deg1] (x) at (0,0) {3};
      \node[deg3] (a) at (1.2,0) {0};
      \node[deg2] (b) at (2.4,0) {1};
      \node[deg1] (z) at (3.6,0) {2};
      \node[deg1] (y) at (1.2,-1.2) {4};
      \draw[edge] (x) -- (a) -- (b) -- (z) (a) -- (y);
    \end{scope}
    \begin{scope}[xshift=5.2cm, yshift=1.0cm]
      \def\cw{0.5}
      \foreach \i/\col in {0/oi-blue, 1/oi-sky, 2/oi-yellow, 3/oi-yellow, 4/oi-yellow}
        \fill[\col!30] (0.1,-\i*\cw-0.22) rectangle (3.95,-\i*\cw+0.22);
      \foreach \w [count=\j from 0] in {1,2,4,8,16}
        \node[font=\scriptsize] at (\j*\cw+0.35,0.55) {$\w$};
      \node[font=\scriptsize] at (3.05,0.55) {$d$};
      \node[font=\scriptsize] at (3.65,0.55) {row};
      \foreach \r [count=\i from 0] in {{0,1,0,1,1},{1,0,1,0,0},{0,1,0,0,0},{1,0,0,0,0},{1,0,0,0,0}} {
        \node[font=\scriptsize, text=black!70] at (-0.2,-\i*\cw) {$\i$};
        \foreach \e [count=\j from 0] in \r
          \node[font=\small] at (\j*\cw+0.35,-\i*\cw) {$\e$};
      }
      \foreach \d/\v [count=\i from 0] in {3/26, 2/5, 1/2, 1/1, 1/1} {
        \node[font=\small] at (3.05,-\i*\cw) {$\d$};
        \node[font=\small] at (3.65,-\i*\cw) {$\v$};
      }
      \draw[thick] (0.15,0.25) -- (0.1,0.25) -- (0.1,-2.25) -- (0.15,-2.25);
      \draw[thick] (2.55,0.25) -- (2.6,0.25) -- (2.6,-2.25) -- (2.55,-2.25);
      \draw[decorate, decoration={brace, amplitude=4pt}, thick]
        (4.05,-0.8) -- (4.05,-2.2)
        node[midway, right=5pt, font=\small] {$2\geq1\geq1$};
    \end{scope}
  \end{tikzpicture}
  \caption{Left: the carbon skeleton of 2-methylbutane with a labelling by
  \(I_5\); dark, medium and light shading mark degree \(3\), \(2\) and \(1\). Right: its adjacency
  matrix, rows and columns in label order; the number above column \(j\)
  is its weight \(2^j\). The columns \(d\) and row give the degree and the row value of each
  vertex. Degrees do not increase down the rows, and the brace shows that the
  row values \(2\geq1\geq1\) of the degree-\(1\) class do not increase either.}
  \label{fig:motivation}
\end{figure}

Let \(G=(V,E)\) be a finite simple undirected graph, and put \(n=|V|\); the
graph need not be connected, and it may have no edges or no vertices. For
\(v\in V\), let \(d(v)\) be its degree and let
\[
  N(v)=\{u\in V:\{u,v\}\in E\}
\]
be its neighbourhood. We write \(u\sim v\) when \(u\) and \(v\) are adjacent.
For each integer \(r\), the degree-\(r\) class is
\(\{v\in V:d(v)=r\}\).

A labelling is a bijection \(f:V\longrightarrow I_n\). A labelling \(f\) is
degree-respecting if it satisfies, for all \(v,w\in V\), the degree condition
\begin{equation}
  d(v)<d(w)\quad\Longrightarrow\quad f(w)<f(v).
\label{eq:degree-order}
\end{equation}

Let \(c:I_n\to\R\) be strictly increasing. Thus
\begin{equation}
  c(0)<c(1)<\cdots<c(n-1).
  \label{eq:strict-weights}
\end{equation}
For a labelling \(f\), write the neighbour sum at \(v\) as
\[
  \sum_{u\in N(v)}c(f(u)).
\]
We say that \(f\) satisfies the equal-degree neighbour-sum condition for \(c\) if, for all
\(v,w\in V\),
\begin{equation}
  d(v)=d(w)\ \text{and}\ f(w)<f(v)
  \quad\Longrightarrow\quad
  \sum_{u\in N(v)}c(f(u))\leq\sum_{u\in N(w)}c(f(u)).
  \label{eq:equal-degree-order}
\end{equation}

\begin{theorem}[Strictly increasing weights]\label{thm:strict-weights}
For every finite simple graph \(G=(V,E)\), with \(n=|V|\), and every strictly
increasing function \(c:I_n\to\R\), there exists a labelling \(f\) satisfying
the degree condition \eqref{eq:degree-order} and the equal-degree neighbour-sum
condition \eqref{eq:equal-degree-order}.
\end{theorem}

\begin{corollary}[Powers of two]\label{cor:powers-of-two}
For every finite simple graph \(G=(V,E)\), with \(n=|V|\), there exists a
labelling \(f\) satisfying the degree condition \eqref{eq:degree-order} and, for
all \(v,w\in V\),
\[
  d(v)=d(w)\ \text{and}\ f(w)<f(v)
  \quad\Longrightarrow\quad
  \sum_{u\in N(v)}2^{f(u)}
  \leq
  \sum_{u\in N(w)}2^{f(u)}.
\]
\end{corollary}

\begin{proof}
Apply \cref{thm:strict-weights} to the strictly increasing function
\(c(i)=2^i\) on \(I_n\).
\end{proof}

\begin{remark}[The two labelling conventions]\label{rem:conventions}
Suppose the alternative convention labels vertices by
\(\{1,\ldots,n\}\) and assigns weight \(2^k\) to label \(k\). For
\(i=k-1\in I_n\), we have
\[
  2^k=2^{i+1}=2\,2^i.
\]
Thus every neighbour sum in the \(\{1,\ldots,n\}\) convention is twice the
corresponding sum in the \(I_n\) convention. The shift \(k=i+1\) preserves the
order of the labels, and the factor \(2\) is common to all vertices, so neither
rule changes.
\end{remark}

\section{Proof of Theorem~\ref{thm:strict-weights}}\label{sec:proof}

Let \(\mathcal A\) be the set of all degree-respecting labellings. The set
\(\mathcal A\) is finite, and it is nonempty: order the vertices by non-increasing degree, choose any order inside each
degree class, and assign the labels in \(I_n\) in that order.

For \(g\in\mathcal A\), define the edge potential
\begin{equation}
  \Phi(g)=\sum_{\{x,y\}\in E}c(g(x))c(g(y)).
  \label{eq:potential}
\end{equation}
The finite nonempty set of real values \(\{\Phi(g):g\in\mathcal A\}\) has a
minimum. Choose \(f\in\mathcal A\) with minimal edge potential.

Take \(v,w\in V\) with \(d(v)=d(w)\) and \(f(w)<f(v)\). Strict increase in
\eqref{eq:strict-weights} gives \(c(f(w))<c(f(v))\). Let \(f'\) be obtained from
\(f\) by exchanging the labels of \(v\) and \(w\). Since
\(d(v)=d(w)\), the exchanged labelling still belongs to \(\mathcal A\).

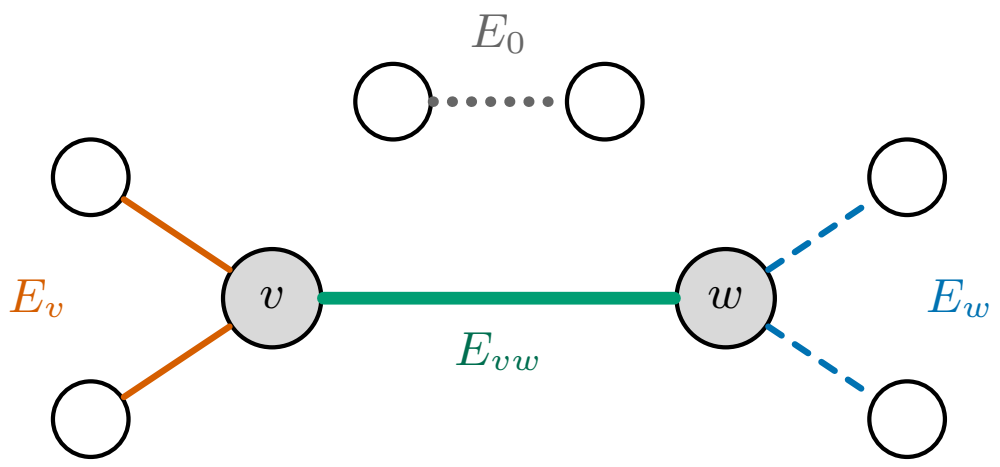
\begin{figure}[ht]
  \centering
  \begin{tikzpicture}[line cap=round, line join=round]
    \node[vtx, fill=black!15] (v) at (0,0) {$v$};
    \node[vtx, fill=black!15] (w) at (3,0) {$w$};
    \node[plain] (v1) at (-1.2,0.8) {};
    \node[plain] (v2) at (-1.2,-0.8) {};
    \node[plain] (w1) at (4.2,0.8) {};
    \node[plain] (w2) at (4.2,-0.8) {};
    \node[plain] (o1) at (0.8,1.3) {};
    \node[plain] (o2) at (2.2,1.3) {};
    \draw[edge, oi-vermillion, very thick] (v) -- (v1) (v) -- (v2);
    \draw[edge, oi-blue, very thick, dashed] (w) -- (w1) (w) -- (w2);
    \draw[edge, oi-green, line width=2.4pt] (v) -- (w);
    \draw[black!60, line width=2pt, dash pattern=on 0pt off 3.5pt]
      (o1) -- (o2);
    \node[font=\small, text=oi-vermillion] at (-1.55,0) {$E_v$};
    \node[font=\small, text=oi-blue] at (4.55,0) {$E_w$};
    \node[font=\small, text=oi-green!70!black] at (1.5,-0.35) {$E_{vw}$};
    \node[font=\small, text=black!60] at (1.5,1.75) {$E_0$};
  \end{tikzpicture}
  \caption{The four edge classes in the comparison of \(\Phi(f')\) with
  \(\Phi(f)\) for the exchanged vertices \(v\) and \(w\): edges incident with
  \(v\) only (\(E_v\), solid), with \(w\) only (\(E_w\), dashed), with both
  (\(E_{vw}\), heavy), and with neither (\(E_0\), dotted).}
  \label{fig:edge-classes}
\end{figure}

To compare \(\Phi(f')\) with \(\Phi(f)\), split the edge set \(E\) into the
edges incident with neither \(v\) nor \(w\) (\(E_0\)), with \(v\) only
(\(E_v\)), with \(w\) only (\(E_w\)), and with both (\(E_{vw}\))
(\cref{fig:edge-classes}). Because the graph is simple,
\(E_{vw}\) is either \(\{\{v,w\}\}\) or empty.

Edges in \(E_0\) have the same contribution before and after the exchange. The
contributions from \(E_v\) and \(E_w\), together with the unchanged contribution
from \(E_{vw}\), give
\begin{equation}
  \begin{aligned}
  &\Phi(f')-\Phi(f)\\
    &\quad={}\bigl(c(f(w))-c(f(v))\bigr)
      \sum_{u\in N(v)\setminus\{w\}}c(f(u))\\
    &\qquad+\bigl(c(f(v))-c(f(w))\bigr)
      \sum_{u\in N(w)\setminus\{v\}}c(f(u))\\
    &\quad=\bigl(c(f(v))-c(f(w))\bigr)
      \left(
        \sum_{u\in N(w)\setminus\{v\}}c(f(u))
        -\sum_{u\in N(v)\setminus\{w\}}c(f(u))
      \right).
  \end{aligned}
  \label{eq:swap-identity}
\end{equation}
Minimality of \(f\) gives \(\Phi(f')-\Phi(f)\geq0\). Since
\(c(f(v))-c(f(w))>0\), \eqref{eq:swap-identity} implies
\begin{equation}
  \sum_{u\in N(v)\setminus\{w\}}c(f(u))
  \leq
  \sum_{u\in N(w)\setminus\{v\}}c(f(u)).
  \label{eq:neighbour-sums-except}
\end{equation}

\paragraph{Adjacent case.}
If \(v\sim w\), then
\[
  \begin{aligned}
  \sum_{u\in N(v)}c(f(u))
    &=\sum_{u\in N(v)\setminus\{w\}}c(f(u))+c(f(w))\\
    &\leq\sum_{u\in N(w)\setminus\{v\}}c(f(u))+c(f(v))\\
    &=\sum_{u\in N(w)}c(f(u)).
  \end{aligned}
\]

\paragraph{Non-adjacent case.}
If \(v\not\sim w\), then
\[
  \sum_{u\in N(v)}c(f(u))
  =\sum_{u\in N(v)\setminus\{w\}}c(f(u))
  \leq\sum_{u\in N(w)\setminus\{v\}}c(f(u))
  =\sum_{u\in N(w)}c(f(u)).
\]
In both cases \(f\) satisfies \eqref{eq:equal-degree-order}, which proves \cref{thm:strict-weights}.

\begin{example}\label{ex:swap}
For the labelling in \cref{fig:motivation}, take \(c(i)=2^i\). The degree
condition \eqref{eq:degree-order} holds, and the only degree class with more
than one vertex is the degree-\(1\) class, with labels \(2,3,4\) and neighbour
sums \(2,1,1\). These sums do not increase with the label, so
\eqref{eq:equal-degree-order} holds as well. The labelling is not the only one:
exchanging the labels \(3\) and \(4\) of the two leaves at the degree-\(3\)
vertex gives a second labelling satisfying both conditions. There are no
others: \eqref{eq:degree-order} forces the labels \(0\) and \(1\) on the
vertices of degree \(3\) and \(2\), and the leaf adjacent to the vertex
labelled \(1\) must carry the label \(2\).

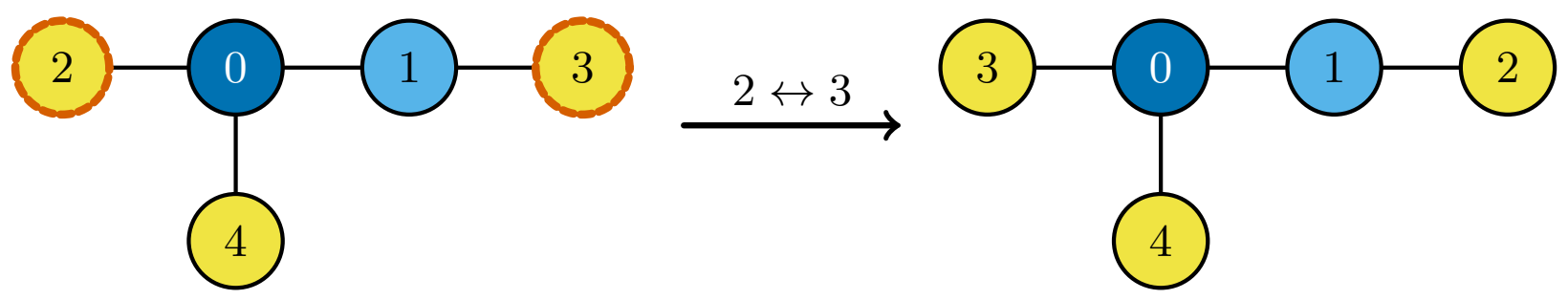
\Cref{fig:swap} starts instead from a degree-respecting labelling that violates
\eqref{eq:equal-degree-order} and shows that exchanging the labels of the
offending pair, the exchange used in the proof above, removes the violation.
\end{example}

\begin{figure}[ht]
  \centering
  \begin{tikzpicture}[line cap=round, line join=round]
    \foreach \dx/\lx/\lz/\sa/\sb/\sc/\tag in
        {0/2/3/1/2/1/before, 6.4/3/2/2/1/1/after} {
      \begin{scope}[xshift=\dx cm]
        \node[deg1] (x-\tag) at (0,0) {\lx};
        \node[deg3] (a-\tag) at (1.2,0) {0};
        \node[deg2] (b-\tag) at (2.4,0) {1};
        \node[deg1] (z-\tag) at (3.6,0) {\lz};
        \node[deg1] (y-\tag) at (1.2,-1.2) {4};
        \draw[edge] (x-\tag) -- (a-\tag) -- (b-\tag) -- (z-\tag)
          (a-\tag) -- (y-\tag);
        \begin{scope}[xshift=1.2cm, yshift=-2.1cm, font=\small]
          \node[anchor=east] at (0,0) {label};
          \node[anchor=east] at (0,-0.5) {sum};
          \foreach \l/\sm [count=\k from 0] in {2/\sa, 3/\sb, 4/\sc} {
            \node at (0.35+\k*0.65,0) {$\l$};
            \node at (0.35+\k*0.65,-0.5) {$\sm$};
          }
          \draw (0.05,-0.25) -- (1.95,-0.25);
        \end{scope}
      \end{scope}
    }
    \node[deg1, draw=oi-vermillion, line width=1.6pt, dash pattern=on 2.2pt off 1.4pt]
      at (x-before) {2};
    \node[deg1, draw=oi-vermillion, line width=1.6pt, dash pattern=on 2.2pt off 1.4pt]
      at (z-before) {3};
    \draw[oi-vermillion, very thick, dashed]
      (1.32,-2.85) rectangle (2.43,-1.85);
    \draw[->, very thick] (4.3,-0.4) -- node[above, font=\small]
      {$2\leftrightarrow3$} (5.8,-0.4);
  \end{tikzpicture}
  \caption{Left: a degree-respecting labelling of the graph of
  \cref{fig:motivation}. The table lists the labels of the
  degree-\(1\) vertices and, below them, their neighbour sums for
  \(c(i)=2^i\); the leaf labelled \(3\) has sum \(2\),
  larger than the sum \(1\) of the leaf labelled \(2\), so
  \eqref{eq:equal-degree-order} fails for the dashed pair. Right: after the two
  labels are exchanged, giving the labelling of \cref{fig:motivation}, the sums
  \(2,1,1\) do not increase with the label, and the edge potential
  \(\Phi\) of \eqref{eq:potential} with \(c(i)=2^i\) drops from \(38\) to \(34\).}
  \label{fig:swap}
\end{figure}

\section{Remarks}\label{sec:remarks}

\Cref{thm:strict-weights,cor:powers-of-two} assert existence, not uniqueness.
A minimiser of the edge potential defined in \eqref{eq:potential} is one labelling among
possibly several satisfying both conditions (\cref{ex:swap}). The two rules do not give a canonical form for graph isomorphism
either: for \(5\leq n\leq7\) some isomorphism classes keep more than one matrix
(\cref{tab:counts}), so the candidates that survive the pruning still have to
be compared. The enumeration of \cite{smooth-pudding-2024} reduces these
comparisons with a variant of Morgan's algorithm \cite{morgan-1965}, which
assigns to each atom of a chemical structure an invariant computed from
iterated sums over neighbours; the surviving candidates are sorted into classes
such that only candidates in the same class can be isomorphic. That filtering
step is separate from the pruning by the two rules studied here.

\Cref{tab:counts} shows the effect of the two rules of the introduction for
\(2\leq n\leq7\). Except in the first and last columns, every entry is a number of adjacency
matrices of size \(n\), that is, of labelled graphs on \(I_n\); an unlabelled
graph \(G\) on \(n\) vertices contributes \(n!/|\mathrm{Aut}(G)|\) distinct
matrices. The columns count, in turn, all such matrices; the skeleton
matrices, whose graph is connected with maximum degree at most \(4\); the
skeleton matrices that also satisfy rule~1, that is, \eqref{eq:degree-order};
and those that satisfy rule~2 as well, that is, \eqref{eq:equal-degree-order}
with \(c(i)=2^i\). The last column counts isomorphism classes of skeletons.

\begin{table}[ht]
  \centering
  \caption{Numbers of adjacency matrices of size \(n\), in total and after each
  successive restriction, and the number of isomorphism classes of connected
  graphs on \(n\) vertices with maximum degree at most \(4\).}
  \begin{tabular}{rrrrrr}
    \hline
    \(n\) & all matrices & skeleton & rule 1 & rules 1 and 2
      & classes\\
    \hline
    2 & 2 & 1 & 1 & 1 & 1\\
    3 & 8 & 4 & 2 & 2 & 2\\
    4 & 64 & 38 & 9 & 6 & 6\\
    5 & 1024 & 728 & 61 & 22 & 21\\
    6 & 32\,768 & 21\,385 & 703 & 100 & 78\\
    7 & 2\,097\,152 & 859\,130 & 10\,544 & 566 & 353\\
    \hline
  \end{tabular}
  \label{tab:counts}
\end{table}

For \(n=7\), rule~1 reduces the \(859\,130\) skeleton matrices to \(10\,544\), by
a factor of about \(81\), and rule~2 reduces these to \(566\), by a further
factor of about \(19\). The \(353\) isomorphism classes then have about \(1.6\)
surviving matrices on average, against about \(2434\) before pruning.

All counts are exact and come from an exhaustive enumeration of the
\(2^{n(n-1)/2}\) matrices. Isomorphism classes are counted by a canonical form,
the lexicographically largest sequence of row values over all relabellings that
keep the degrees non-increasing; two matrices are isomorphic exactly when their
canonical forms coincide. Applied to every skeleton
matrix, it confirms that each class keeps at least one matrix under both rules, as
\cref{cor:powers-of-two} guarantees.
The class counts agree with those reported in \cite{smooth-pudding-2024}.

\Cref{thm:strict-weights,cor:powers-of-two}, with their definitions and proofs,
have been formalised and mechanically checked in Lean~4, a proof assistant
\cite{degree-labeling-lean-2026}.
The Lean development represents the finite vertex set as a finite type (a type
with finitely many elements), represents a labelling as a bijection to
\(I_n\), and represents the edge potential by a finite sum over the edge set.
Its proof dependencies are three foundational principles: propositional extensionality,
which identifies logically equivalent propositions; choice, which selects an
element from a nonempty collection; and quotient soundness, which makes
related elements equal in the quotient. Their Lean names are
\texttt{propext}, \texttt{Classical.choice}, and \texttt{Quot.sound}; these are
dependencies of the formal proof, not hypotheses on the graph or on the weights.

\bibliographystyle{amsplain}
\bibliography{../references}

\end{CJK*}
\end{document}